\documentclass{notices}
\usepackage{amsfonts,amssymb,amsmath,amscd}
\usepackage{url}
\usepackage{graphicx,epstopdf,color}
\usepackage{amsrefs}

\newcommand{\R}{\mathbb{R}}

\renewcommand{\epsilon}{\varepsilon}

\renewcommand{\le}{\leqslant}

\renewcommand{\ge}{\geqslant}

\title{Towards a long-range theory of capillarity}

\author{Enrico Valdinoci\affil{Department of Mathematics and Statistics,
University of Western Australia, 35 Stirling Highway,
WA6009 Crawley, Australia}}

\begin{document}
\maketitle

{\sl
We recall the classical theory of capillarity, describing the shape of a liquid droplet in a container, and present a recent
approach which aims at accounting for long-range particle interactions.

This nonlocal setting recovers the classical notion of surface tension in the limit. We provide some regularity results and the determination of the contact angle, supplied with various asymptotics. }

\section{The classical theory of capillarity}

\subsection{The surface tension as an average of long-range molecular forces}
The classical capillarity theory aims at understanding the displacement of a liquid droplet in a container in view of
surface tension.

In spite of this elementary formulation, capillarity
is a very delicate issue and its study emerges in a number of fields, including the calculus of variations, geometric measure theory, partial differential equations, mathematical physics, material sciences, biology and chemistry. A complete understanding of the intriguing patterns related to capillarity will require the virtuous blend of different ideas coming both from mathematics and from the applied sciences.

Indeed, surface tension is a complex phenomenon arising as the average outcome of the attractive forces between molecules (such as cohesion and adhesion)
and accounts for the interfaces between the droplet, the air, and the container.

The classical capillarity theory models the surface tension as a local average of
intermolecular forces which in principle possess long-range contributions. Namely,
it is classically assumed that away from the interfaces between the different materials
the cohesive forces average out (because each molecule is pulled pretty much equally in every direction by the other molecules) and accordingly the surface tension is seen as a force
concentrated at the interface.

\begin{figure}[h]
    \centering
    \includegraphics[width=6cm]{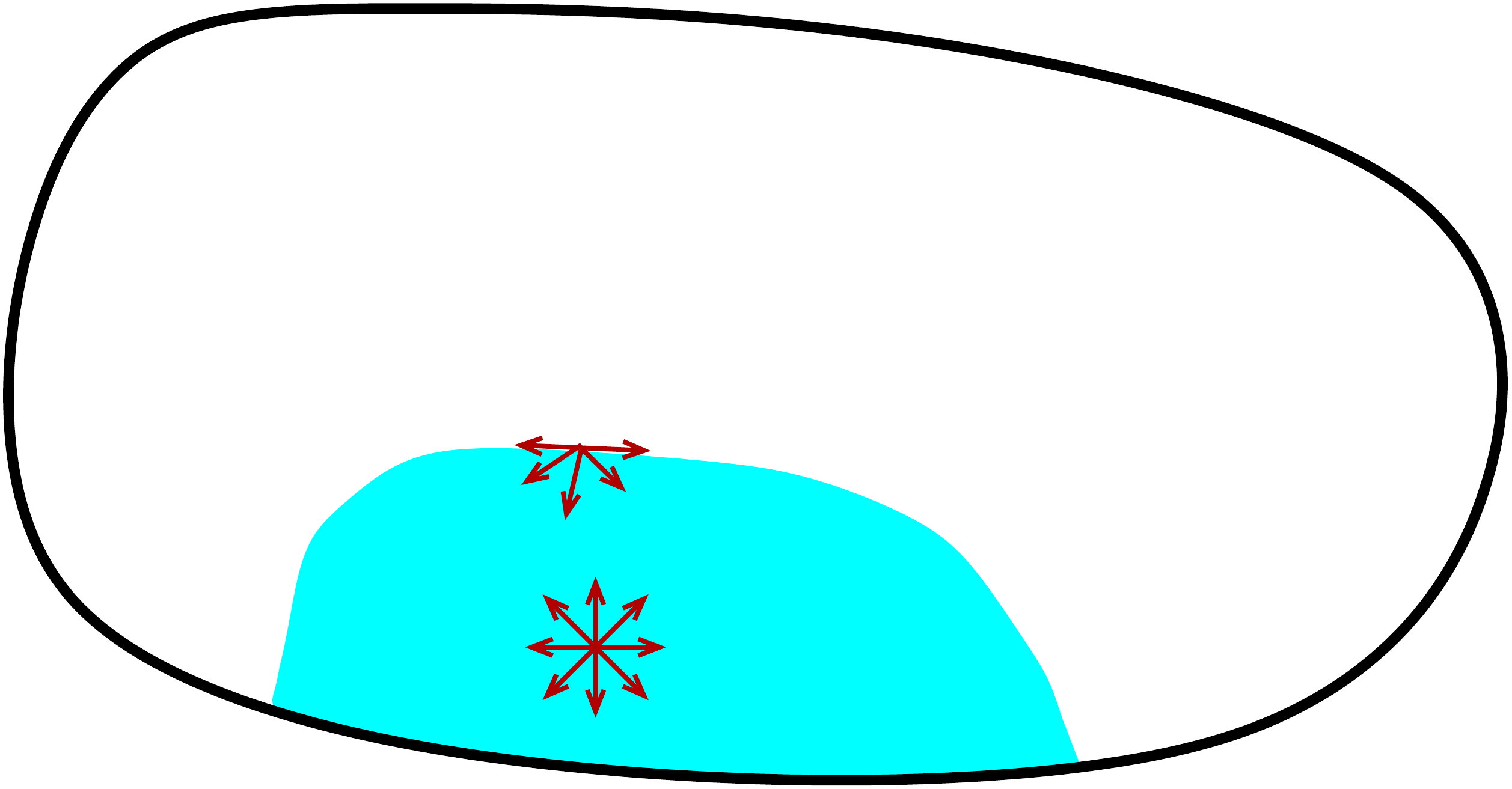}
    \caption{Averaging out long-range interactions to approximate the net force by a surface tension.}
    \label{GOC1}
    \end{figure}

The net force on the interface is not zero, since, for a smooth separation of media, focusing on the ``infinitesimal forces'' around a given point, one expects the cohesion forces to act on ``one side'' of the interface (or, more precisely, that the the forces acting on one side of the interface are different than the ones on the other side,
see Figure~\ref{GOC1}). In view of these considerations, the interface in the classical theory is expected to shrink into the minimum surface area possible.

\begin{figure}[h]
    \centering
    \includegraphics[width=6cm]{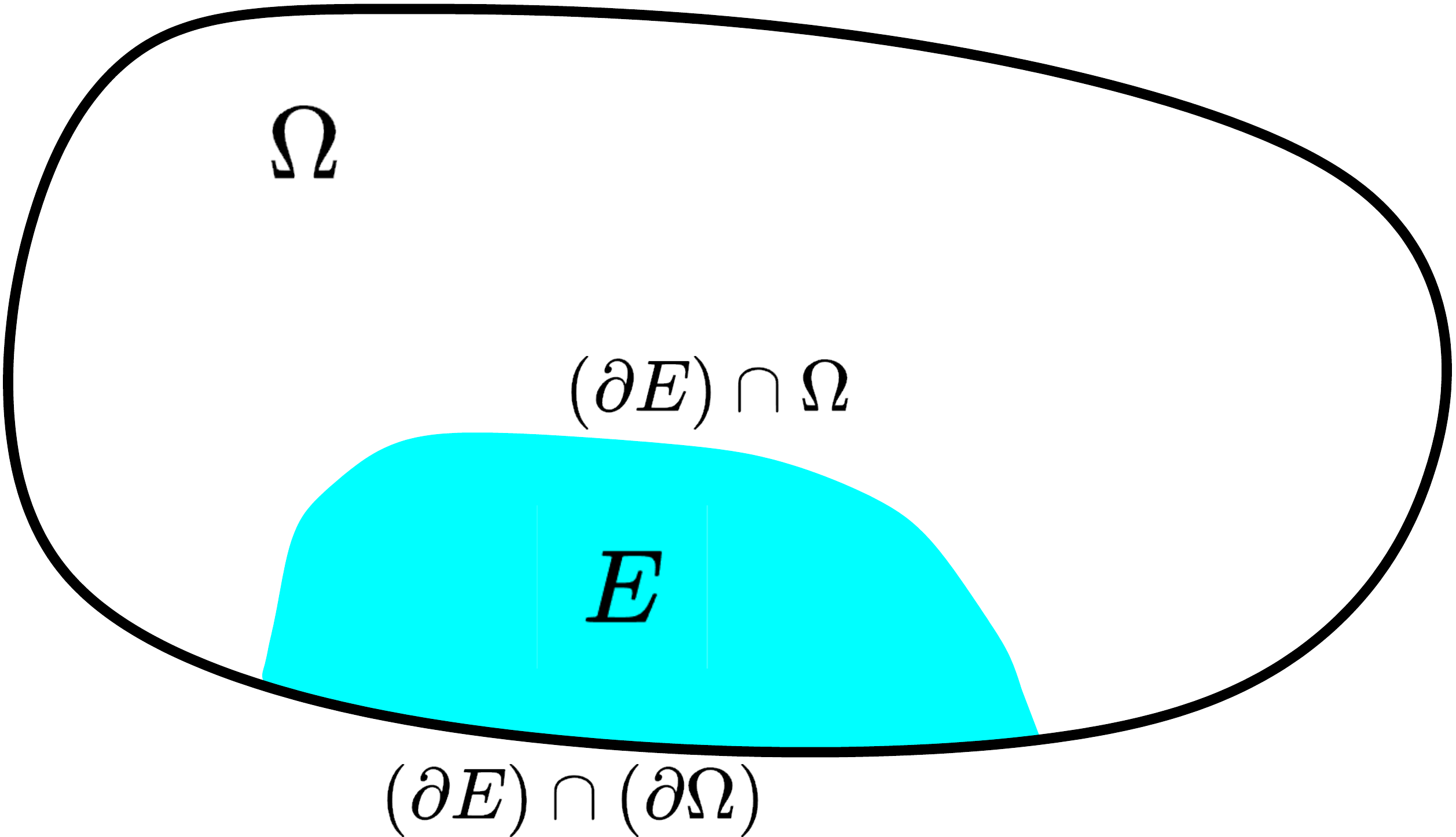}
    \caption{Different contributions to the classical capillarity energy functional.}
    \label{GOC2}
    \end{figure}

{F}rom the mathematical standpoint, this analysis suggests that
the formation of droplets in a container~$\Omega\subset\R^n$ (supposed to be open and smooth) is given
by the minimization of a ``perimeter-like'' energy functional.
Specifically, a droplet is modeled by a set~$E\subset\Omega$, with a prescribed volume~$m$, and its capillarity energy (up to normalizing constants) is given by the functional
\begin{equation}\label{PAKSMDP} {\mathcal{E}}_\sigma(E):=|(\partial E)\cap\Omega|+\sigma \,|(\partial E)\cap(\partial \Omega)|.\end{equation}
One can also add to the functional external forces, such as gravity, but for the sake of simplicity we will focus here on the simplest possible scenario. Here above,
the quantity~$|(\partial E)\cap\Omega|$ stands for the $(n-1)$-dimensional surface area
of the set~$(\partial E)\cap\Omega$ (which is the surface of the droplet inside the container)
and~$|(\partial E)\cap(\partial \Omega)|$ for the $(n-1)$-dimensional surface area
of the set~$(\partial E)\cap(\partial \Omega)$ (which is the wet region of the container, see Figure~\ref{GOC2}).

The parameter~$\sigma\in\R$ is called ``relative adhesion coefficient'' (for a reason that will be clarified in~\eqref{SACSAS} below) and it ``weighs'' the importance of the intermolecular forces among the liquid particles versus the ones between the liquid and the container.

It is thereby natural to consider the minimizers of the functional in~\eqref{PAKSMDP} as the equilibrium configurations of a droplet~$E$ in a container~$\Omega$.
    
\subsection{The relative adhesion coefficient}
It is natural to restrict the values of~$\sigma$ in~\eqref{PAKSMDP} to the interval~$[-1,1]$.
For instance, let us show that if there exists a minimizer containing a ball centered at the boundary of the container and whose complement inside~$\Omega$ contains a ball, then necessarily
\begin{equation}\label{IS-1}
\sigma\le1
.\end{equation}

\begin{figure}[h]
    \centering
    \includegraphics[width=6cm]{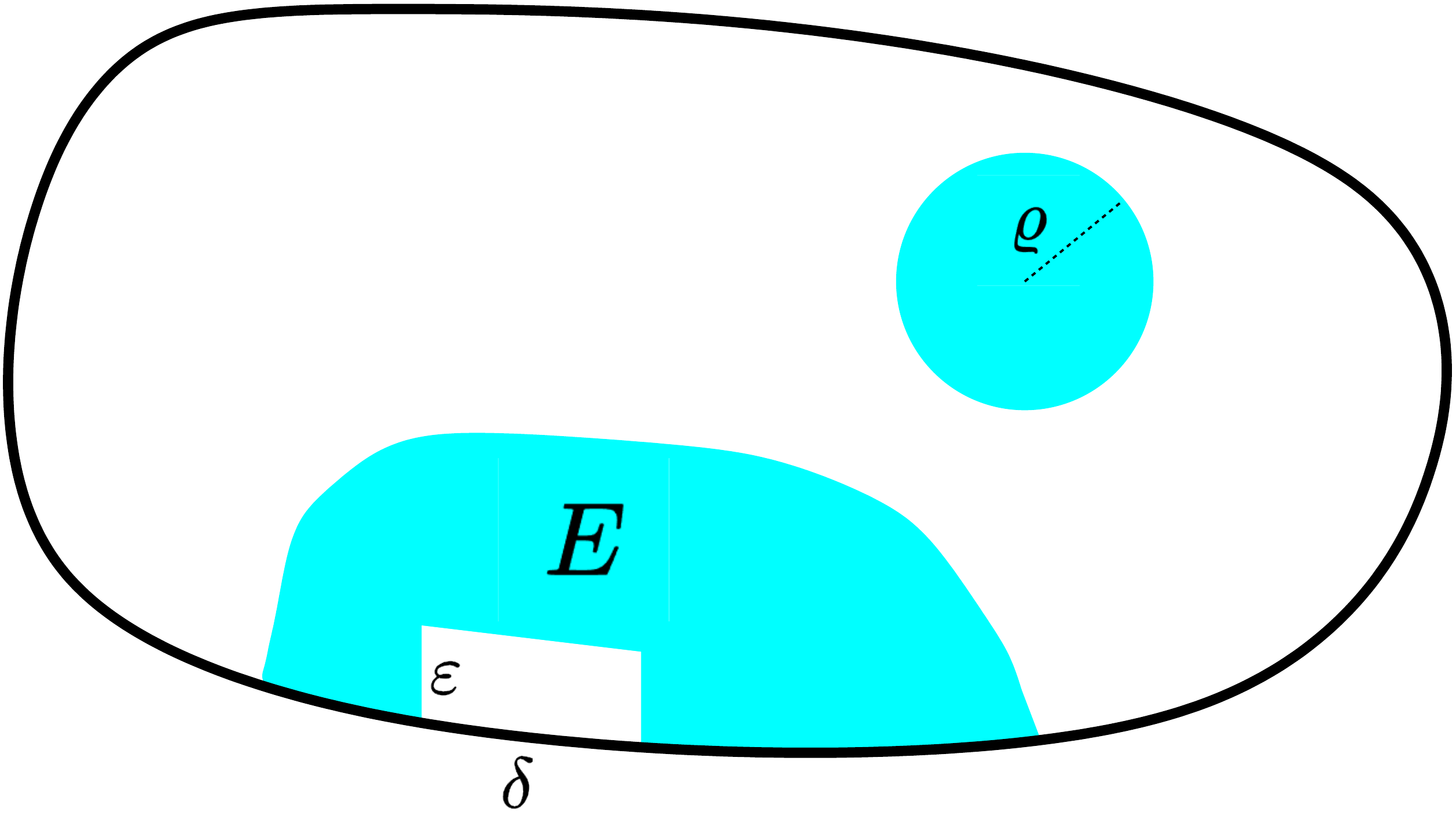}
    \caption{Constructing a competitor for the minimizer.}
    \label{GOC3}
    \end{figure}

To check this, one can create a competitor by digging out some mass near the boundary of~$\Omega$,
as shown in Figure~\ref{GOC3}, in the shape of an ``infinitesimal'' rectangle which is a small deformation of~$(0,\delta)^{n-1}\times(0,\varepsilon)$, and (to maintain the total mass of the droplet) by adding an ``infinitesimal'' ball of volume~$\varepsilon\delta^{n-1}$
(hence, of surface area~$c_n\,(\varepsilon\delta^{n-1})^{\frac{n-1}n}$ for some~$c_n>0$).

The minimality condition gives that 
\begin{eqnarray*} \sigma\delta^{n-1} \le \delta^{n-1}+2(n-1)\varepsilon\delta^{n-2}+c_n\,(\epsilon\delta^{n-1})^{\frac{n-1}n}
\end{eqnarray*}
Dividing by~$\delta^{n-1}$ and sending~$\varepsilon\searrow0$, we find that
\begin{eqnarray*} \sigma \le\lim_{\varepsilon\searrow0} 1+2(n-1)\varepsilon\delta^{-1}+c_n\,\varepsilon^{\frac{n-1}n} \delta^{\frac{1-n}n}=1,
\end{eqnarray*}
proving~\eqref{IS-1}. 

After this, noticing that~${\mathcal{E}}_\sigma(E)=|\partial\Omega|-{\mathcal{E}}_{-\sigma}(E^c)$, where~$E^c:=\Omega\setminus E$, one also infers that it is natural to suppose that $$\sigma\ge-1.$$

\begin{figure}[h]
    \centering
    \includegraphics[width=4.25cm]{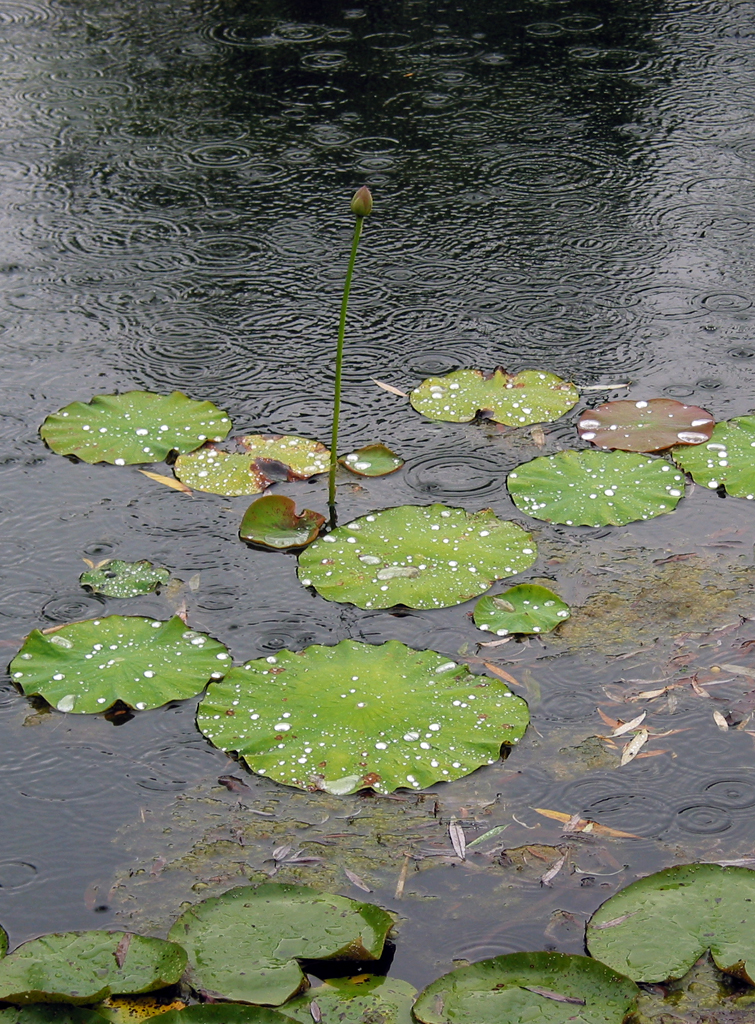}
    \caption{Lotus leaves on pond in rain. Photo by ArchiKat, taken in Erholungspark Marzahn, Berlin,  licensed under the Creative Commons Attribution 3.0 Unported.}
    \label{GOC4}
    \end{figure}

When~$\sigma=1$, we have that the capillarity functional in~\eqref{PAKSMDP} reduces to the surface area of~$E$ and consequently, by the isoperimetric inequality, the minimum for small enough volume of the droplet
is a ball (actually, any ball with the prescribed volume~$m$, placed wherever in~$\Omega$).
This situation is called ``completely non-wetting'', since the 
droplet detaches from the
boundary of the container, which thus remains
dry (up to at most a point) and corresponds to the so-called\footnote{Capillarity
is a very complex phenomenon, also depending on the high heterogeneity of the surfaces
and on their specific microscopic scales (this is also the reason for which the study of capillarity is deeply intertwined with that of nanomaterials). The particular
hierarchies of nanostructures involved produce, for example, a different
behavior between the lotus leaf (in which the droplet rolls off easily)
and the petals of the rose (whose
micropapillae and nanofolds retain the droplet, which typically
cannot roll off even if the petal is turned upside down). This phenomenon is called ``petal effect''. Here we will not dive into the analysis of the nanostructures
that are at the basis of the capillarity theory from the point of view of materials science,
but rather look at mathematically simplified models to deal with molecule interactions.}
``lotus effect'' (see Figure~\ref{GOC4}).

This liquid repellency phenomenon has important biological consequences, since it provides
a self-cleaning mechanism for the leaves of several plants, in which dirt particles are picked up and washed out by water droplets, as well as self-drying procedures for the wings of several insects.
The lotus effect is also used in a number of practical applications, including bathers, rashguards, windshields, paints, antibacterial surfaces, etc. 

An interesting application of hydrophobic materials also consists in enhancing the droplet rebounds on a surface, reducing the viscous dissipation in view of a non-wetting situation, see e.g.~\cite[Figure~10]{QUERE200232}.

Instead, when~$\sigma=-1$, the functional in~\eqref{PAKSMDP} boils down to~$|(\partial E)\cap\Omega|-|(\partial E)\cap(\partial \Omega)|$,
which has the same minima as~$|(\partial E^c)\cap\Omega|+|(\partial E^c)\cap(\partial \Omega)|$. Therefore, if the volume of~$E$ is sufficiently large (hence the volume of~$E^c$ is sufficiently small), the miminum is achieved when~$|(\partial E^c)\cap(\partial \Omega)|=0$
and, by isoperimetric inequality, when~$E^c$ is a ball (hence, $E$ is the complement of a ball).
This is a ``perfectly wetting'' scenario, where the liquid
completely sticks to the container, which also has a number of practical applications,
such as lubricants (e.g., contact lenses, which need to remain wet by tear fluid not to hurt the eyes) and anti-fogging swimming goggles (in which 
the humidity built in the goggles do not form droplets on
the inner surface, but rather spreads on it as a thin, transparent film).\medskip

In general, the sign of~$\sigma$ distinguishes between ``hydrophobic'' materials,
in which~$\sigma>0$, thus producing a ``repulsive'' effect
between the container and the droplet as an outcome of the minimization
of the energy functional in~\eqref{PAKSMDP},  and ``hydrophilic'' materials,
in which~$\sigma<0$, corresponding to an ``attractive'' effect
induced by the negative sign in~\eqref{PAKSMDP}.
As we will see below on page~\pageref{ONPA}, this distinction also pops up when considering the shape of the droplet,
and specifically the angle formed between the droplet and the container.

\subsection{Minimality conditions}

Minimizers (and, more generally, critical points) of the energy functional in~\eqref{PAKSMDP}
have\footnote{At this level, we are implicitly assuming minimizers to be ``smooth enough'' for the necessary computations to make sense. See e.g.~\cite[Chapter~19]{MR2976521}
and the reference therein for a more precise description of the classical capillarity theory.} constant mean curvature, since 
one can consider perturbations which preserve volume and keep~$(\partial E)\cap(\partial\Omega)$
as it is (notice that these perturbations yield that the energy term~$|(\partial E)\cap\Omega|$
must be locally stationary with respect to Lagrange multipliers induced by the volume constraint).

Interestingly, energy perturbations involving the contact point
produce a second necessary condition enjoyed by minimizers (and, more generally, critical points) in terms of the so-called ``contact angle'', which is defined as the angle between the tangent to the droplet and that to the container at the contact point.
Remarkably, this contact angle, which will be denoted by~$\vartheta\in[0,\pi]$, possesses a simple relation with the relative adhesion coefficient~$\sigma$, namely
\begin{equation}\label{YOU} \cos\vartheta=-\sigma.\end{equation}
This relation is known as ``Young's Law'', named after
polymath\footnote{Besides discovering the law of the contact angle, Young made several
notable contributions to a number of topics,
including the decipherment of the demotic script, a method of tuning musical instruments,
and several studies in linguistics (introducing the term Indo-European languages),
see~\cite{Y}.} Thomas Young.

\begin{figure}[h]
    \centering
    \includegraphics[width=6cm]{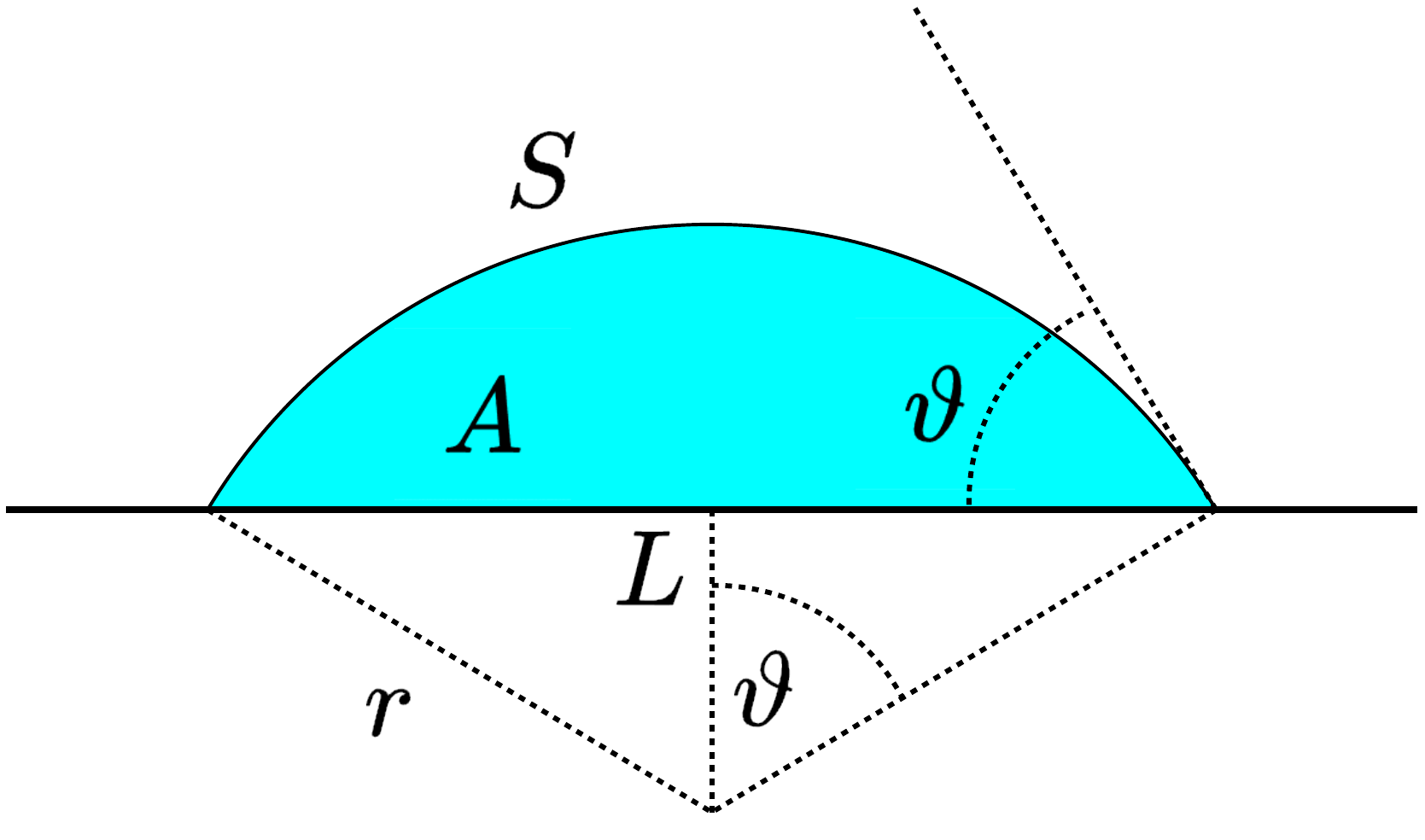}
    \caption{A spherical cap to justify Young's Law.}
    \label{GOC41}
    \end{figure}

To justify, at least heuristically, Young's Law in~\eqref{YOU}, we consider the planar case~$n=2$ and take an ``infinitesimal'' round droplet of radius~$r$.
At this small scale, we identify~$\Omega$ with a halfspace and we recall that the
length of the circular segment~$S$ is~${2\vartheta }r$,
the area of the circular segment~$A$ is~$r^2 \left(\vartheta -\frac{\sin (2\vartheta)}2 \right)$ and the chord length~$L$ of the circular segment is~$2r\sin\vartheta$, see Figure~\ref{GOC41}.

By the volume constraint, we have that $$ m=r^2 \left(\vartheta -\frac{\sin (2\vartheta)}2 \right)$$
and therefore
$$ r=\sqrt{\frac{2m}{2\vartheta -\sin (2\vartheta)}}.$$
In this way, the capillarity energy functional in~\eqref{PAKSMDP} can be written, as a function of~$\vartheta$, in the form
\begin{eqnarray*}&&
S+\sigma L=2\vartheta r+2\sigma r\sin\vartheta\\&&\quad
=\left(2\vartheta +2\sigma \sin\vartheta\right)
\sqrt{\frac{2m}{2\vartheta -\sin (2\vartheta)}}\\&&\quad
=2\left(\vartheta +\sigma \sin\vartheta \right)\sqrt{\frac{m}{\vartheta -\sin \vartheta\cos\vartheta}},
\end{eqnarray*}
whose derivative vanishes if and only if~${\sigma + \cos \vartheta}=0$, in agreement with Young's Law~\eqref{YOU}.\medskip

\begin{figure}[h]
    \centering
    \includegraphics[width=6cm]{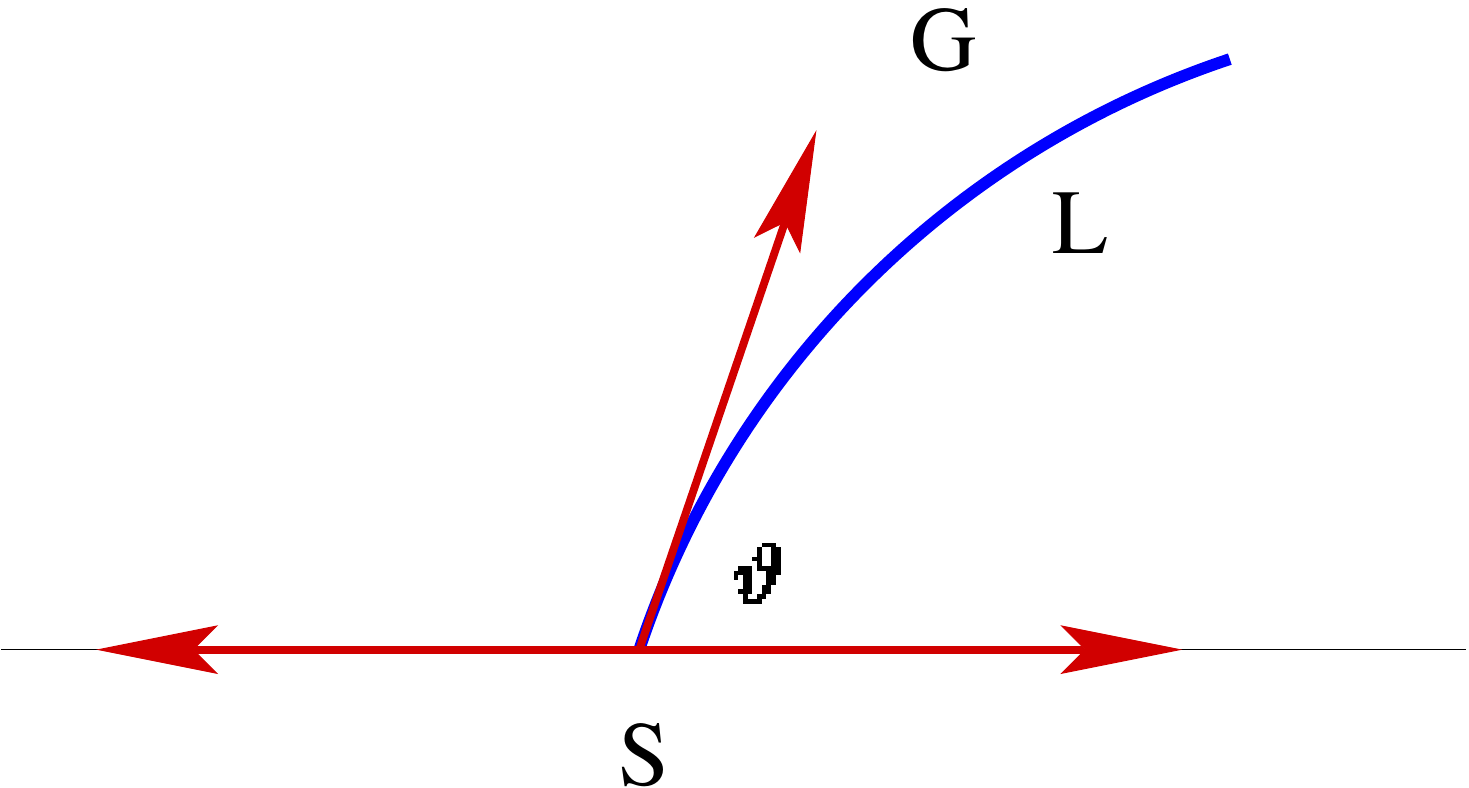}
    \caption{Projecting tension forces at a contact point.}
    \label{GOC8}
    \end{figure}

It is worth observing that, in light of Young's Law,
the role of~$\sigma$ as a relative adhesion coefficient also emerges from simple
physical considerations related to force balance.
Specifically, if we denote by~$\gamma_{{\text{SG}}}$ the interfacial tension between the solid and gas, by~$\gamma_{{\text{SL}}}$ the interfacial tension between the solid and liquid,
and by~$\gamma_{{\text{LG}}}$ the interfacial tension between the liquid and gas, at a contact point between the droplet and the container the balance of forces (see Figure~\ref{GOC8})
should give
$$ \gamma_{\text{SG}}=\gamma_{\text{SL}}+\gamma_{\text{LG}}\cos \vartheta $$
and therefore, by~\eqref{YOU},
\begin{equation}\label{SACSAS} \sigma=-\cos \vartheta=\frac{\gamma_{\text{SL}}-\gamma_{\text{SG}}}{\gamma_{\text{LG}}},\end{equation}
which highlights the role of~$\sigma$ as a ratio of different tension forces.\medskip

One of the consequences of Young's Law is that the
relative adhesion coefficient~$\sigma$ clearly relates the geometry of the droplet
with the water attraction of the container. Indeed, by~\eqref{YOU}, surfaces with positive
relative adhesion coefficient corresponding to hydrophobic materials
produce droplets with a contact angle~$\vartheta$ lying in~$\left(\frac\pi2,\pi\right]$,
while surfaces with negative
relative adhesion coefficient corresponding to hydrophilic materials
produce droplets with a contact angle~$\vartheta$ lying in~$\left[0,\frac\pi2\right)$
(hence, the ``flatter'' the droplet, the highest the attraction to the container, \label{ONPA}
for instance, water droplet on a hydrophobic lotus leaf have been showing contact angles of about~$148^\circ$, while hydrophilic coatings are available on the market warranting a
contact angle of about~$10^\circ$).

\begin{figure}[h]
    \centering
    \includegraphics[width=6cm]{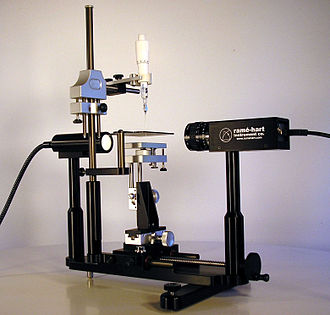}
    \caption{Image of a contact angle goniometer. Photo by Ramehart,  licensed under the Creative Commons Attribution 3.0 Unported.}
    \label{GOC17}
    \end{figure}

Given its topical physical relevance as an indication of the wettability of the container,
the experimental measure of the contact angle utilizes nowadays sophisticated devices,
named ``contact angle goniometers'', which employ high resolution cameras, see Figure~\ref{GOC17}. In practice, especially in the presence of lubricants,
the contact angle of a droplet may undergo swift changes in the vicinity of the surface:
for example, lubricant skirts can produce an ``apparent contact angle'', also calling for a reformulation of the capillarity model under consideration, see e.g.~\cite[Figure~2]{skir}.

\section{A nonlocal capillarity theory}

\subsection{The surface tension as a pointwise interaction} In view of the state of the art on the classical theory of capillarity and inspired by the intense development experienced by nonlocal equations in the recent past, we aim at presenting a new theory of capillarity in which the classical, local notion of surface tension was replaced by a nonlocal energy, accounting for long-range particle interactions. Though, in principle, one
wishes to have a unified description of a physical phenomenon accounting for all scales
(e.g., from the macroscopic description of interfacial regions of scale-invariant surface type
to the nanoscale in which thin films become comparable to the size of the molecules),
the full description of all intermolecular forces in this setting could be extremely challenging, possibly unpractical, and it is a common pragmatic principle to often rely on simplified phenomenological models.
Hence, as an initial model case, we suppose the interactions to be long-range, but invariant under translations and rotations.

Specifically, inspired by the theory of nonlocal minimal surfaces put forth in~\cite{MR2675483}, given two (measurable) disjoint sets~$X$, $Y\subseteq\R^n$
and an interaction kernel~$K:\R^n\times\R^n\to[0,+\infty]$,
we considered in~\cite{MR3717439}
the point-set interaction
$$ {\mathcal{I}}(X,Y):={\iint_{X\times Y}K(x-y)\,dx\,dy }.$$

The prototypical example of kernel we have in mind was the positively homogeneous one, given by
\begin{equation}\label{KA2s} K_s(x)={\frac{s\,(1-s)}{|x|^{n+s}}},\end{equation}
for some~$s\in(0,1)$ (here, the factor~$s\,(1-s)$ is taken as a normalization to ``stabilize''
the limits as~$s\nearrow1$ and~$s\searrow0$).

If not otherwise specified, we will consider kernel that are (sufficiently smooth outside the origin and) comparable with~$K_s$, i.e. we suppose that~$K\in C^1(\R^n\setminus\{0\})$ and there exist~$\varrho_0>0$ and~$\Lambda\ge1$ such that, for all~$x\in\R^n\setminus\{0\}$,
\begin{equation}\label{ILK} \begin{split}&\frac{\chi_{B_{\varrho_0}}{(x)}}{\Lambda\,|x|^{n+s}}\le K(x)\le
\frac{\Lambda}{|x|^{n+s}}\\
{\mbox{and}}\qquad&|\nabla
K(x)|\le\frac{\Lambda}{|x|^{n+s+1}}.\end{split}\end{equation}

This setting allows one to consider long-range interactions between points of the droplet~$E$ with points outside the droplet but inside the container~$\Omega\setminus E$:
namely, the local term~$|(\partial E)\cap\Omega|$ in~\eqref{PAKSMDP} describing the surface tension between the liquid~$E$ and the air in the container~$\Omega\setminus E$ is now replaced by a nonlocal term of the form~${\mathcal{I}}(E,\Omega\setminus E)$.

Similarly, the local term~$|(\partial E)\cap(\partial\Omega)|$ in~\eqref{PAKSMDP} describing the surface tension between the liquid~$E$ and the external of the container~$\R^n\setminus\Omega$ is now replaced by a nonlocal term of the form~${\mathcal{I}}(E,\R^n\setminus\Omega)$.

In this sense, the nonlocal functional
\begin{equation}\label{PAKSMDPs} {\mathcal{E}}_{\sigma,K}(E):=
{\mathcal{I}}(E,\Omega\setminus E)+\sigma\,
{\mathcal{I}}(E,\R^n\setminus\Omega)\end{equation}
can be considered as a long-range counterpart of the capillarity energy functional in~\eqref{PAKSMDP} 
with the first term mimicking the molecular interactions of the fluid with the gas, the second with the solid (again, external forces can be accounted for as well, but we focus here on the simplest possible scenario).
In the spirit of~\cite{MR3586796}, we have that, for the prototypical interaction kernel in~\eqref{KA2s}, as~$s\nearrow1$, the nonlocal functional~${\mathcal{E}}_{\sigma,K}$ in~\eqref{PAKSMDPs} recovers the classical capillarity energy functional~${\mathcal{E}}_{\sigma}$ in~\eqref{PAKSMDP}, up to normalizing constants
(this is physically important, since one of the desirable features of long-range alternative
to the classical surface tension is that they recover comparable macroscopic properties).

The nonlocal capillarity functional~${\mathcal{E}}_{\sigma,K}$ presents certainly a number of new technical and conceptual difficulties.
First of all, we can certainly expect that the classical ``local calculus'', based on computing derivatives, becomes much less effective that in the classical case
(roughly speaking, integration is generally harder than differentiation!)
and we can expect not to be able to write many ``explicit'' solutions of the problems under consideration.

Furthermore, the ``cut-and-paste methods'', such as the ones sketched in Figures~\ref{GOC2} and~\ref{GOC3},
in which one moves mass around to build a convenient competitor become trickier in the nonlocal case, since one can expect that bespoke computations will be needed to
carefully account for all the energy contributions.

In particular, given the possible singularities and the decay properties of the interaction kernel under consideration, one can expect that suitable integral cancellations
have to be attentively spotted in order
to ``average out'' the ``microscopic'' interactions which have little impact on a ``macroscopic'' effect, and a deep understanding of the notion of ``effective scale'' may be required for this kind of analysis (also in consideration of the ``fat tails'' exhibited at infinity by kernels
with polynomial decay such as the one in~\eqref{KA2s}).

These difficulties are compensated however by several advantages. 
The use of nonlocal methods, for instance, may reduce approximation errors and improve numerical stability, see e.g.~\cite[pages~170--171]{MR3824212}. 

\subsection{The structure of the minimizers of the nonlocal capillarity functional}

Given the prominent role of the mean curvature in the description of classical capillarity minimizers, we can expect that, to analyze the geometric structure of the minimizers
of the nonlocal capillarity functional,
it comes in handy to consider the notion of nonlocal mean curvature of a set~$ E$ at the point~$x\in\partial E$ with respect to the kernel~$K$, namely
\begin{equation}\label{NMC} {\mathcal{H}}^K_{ E}(x):=\int_{\R^n}\Big(\chi_{\R^n\setminus E}(y)-\chi_E(y)\Big)\,K(x-y)\,dy.\end{equation}
We remark that this integral converges in the principal value sense as soon as E is~$C^{1,\alpha}$ near~$x$
with~$\alpha\in(s,1)$. This suggests that it is worth considering the ``regular part
of~$\partial E$'', namely the collection of points~$x\in\overline{(\partial E)\cap\Omega}$ such that~$\partial E$ is
a $C^{1,\alpha}$-hypersurface with boundary near~$x$
with~$\alpha\in(s,1)$ and whose boundary points are in~$\partial\Omega$,
for which the nonlocal mean curvature in~\eqref{NMC} is well-defined in the classical sense
(for simplicity, we present the main results assuming without mention this regularity assumption, but more general versions are possible too).

To the casual eye, one could develop the nonlocal theory of capillarity along the lines of the classical case, by replacing ``mean curvature'' by ``nonlocal mean curvature'' and leaving the rest unchanged, but we will see the situation is slightly subtler and, in fact,
more complicated geometries may arise in the nonlocal setting, due to remote interactions. Particularly, differently from the classical case, nonlocal capillarity minimizers are
not portions of balls, not even in dimension~$2$, not even when the container is a halfplane, not even when~$\sigma=0$, see~\cite[Appendix~1]{MR3707346}.

To appreciate the shape of the minimizers, 
one can look at the Euler-Lagrange equation associated with the nonlocal capillarity energy functional (see Theorem~1.3 in~\cite{MR3717439}), yielding that, for
a critical set~$E$ of the nonlocal capillarity functional in~\eqref{PAKSMDPs} and~$
x\in\partial E$,
\begin{equation}\label{el} {\mathcal{H}}^K_{ E}(x)+(\sigma-1)\int_{\R^n\setminus\Omega}K(x-y)\,dy={\rm const}.\end{equation}
Interestingly, this Euler-Lagrange is structurally different from its local counterpart, since it combines a geometric prescription induced by the nonlocal mean curvature with an interaction term coming from the exterior of the container and depending on the relative adhesion coefficient~$\sigma$ (unless~$\sigma=1$, in which case the nonlocal capillarity functional reduces to the nonlocal perimeter in~\cite{MR2675483} and, coherently with this, the corresponding Euler-Lagrange equation boils down to a prescribed mean curvature equation).
Also, when the interaction kernel is the prototypical one in~\eqref{KA2s}, the
additional interaction term in~\eqref{el} vanishes as~$s\nearrow1$, thus recovering the constant mean curvature condition in the limit.
Hence, on the one hand,
differently from the classical case, the relative adhesion coefficient~$\sigma$ appears in the nonlocal Euler-Lagrange equation in~\eqref{el}, but, on the other hand, this dependence vanishes as~$s\nearrow1$.
\medskip

The interior regularity of the minimizers of the nonlocal capillarity functional can be
obtained in the light of the theory of almost minimizers. In particular, exploiting also previous results in~\cites{MR3090533, MR3107529, MR3331523} one obtains smoothness up to a possible singular set of Hausdorff dimension at most~$ n -3$ (and~$n-8$ when the fractional parameter is sufficiently close to~$1$; it is still open to determine whether this is sharp
and to construct any example of singular set). In this way (see Theorem~1.6 in~\cite{MR3717439}),
if~$E$ is a minimizer for the nonlocal capillarity functional in~\eqref{PAKSMDPs}
with the interaction kernel~$K$ is as in~\eqref{KA2s},
then the regular part of~$\partial E$ is
a smooth hypersurface and the singular set has
Hausdorff dimension at most~$n - 3$ (in particular, $\partial E$ is smooth in dimension~$2$).

Moreover, there exists~$\epsilon_0\in(0,1)$ such that if~$s\in(1-\epsilon_0,1)$, 
then the singular set has
Hausdorff dimension at most~$n - 8$ (in particular, in this case~$\partial E$ is smooth in dimension up to~$7$).

\subsection{Nonlocal Young's Law} Having analyzed the basic interior properties of the nonlocal capillarity minimizers, it is also interesting to understand their behavior at the boundary, and, especially, to understand the nonlocal counterpart of the contact angle prescription.

To this end, we first point out that, while a complete boundary regularity theory has not been established yet, nonlocal minimizers satisfy energy and density estimates uniformly up to the boundary. More precisely (see Theorem~1.7 in~\cite{MR3717439}),
if~$E$ is a minimizer for the nonlocal capillarity functional in~\eqref{PAKSMDPs}
with interaction kernel~$K$ is as in~\eqref{KA2s},
then, for small~$r>0$, for all~$p\in\R^n$ we have that
$$ {\mathcal{I}}\big(E\cap B_r(p),\,\R^n\setminus(E\cap B_r(p))\big)\le Cr^{n-s}$$
and, for all~$p\in\overline{(\partial E)\cap\Omega}$,
$$ c\le\frac{|E\cap B_r|}{|E_r|}\le C,$$
with constants~$C>c>0$ depending only on~$n$, $s$, $\sigma$ and~$\Omega$.

To rigorously speak about a contact angle, it is also convenient to perform larger and larger dilations and ``zoom in'' at a point of~$(\partial\Omega)\cap(\partial E)$.
This is a delicate issue, since the previous results do not exclude, in principle,
the formation of small fractal microstructures of~$\partial E$ in the vicinity of~$\partial \Omega$ which would prevent us to speak about a well-defined contact angle.
To circumvent this difficulty, one notices that, for~$K$ as in~\eqref{KA2s}
and~$\Omega$ a halfspace,
the nonlocal capillarity functional 
admits a suitable ``extension problem'' (see~\cite[Proposition~1.1]{MR4404780}), which in turn possesses a convenient boundary monotonicity formula
(see~\cite[Theorem~1.2]{MR4404780}).

This strategy leads to the existence, up to subsequences, of blow-up limits
which are cones (i.e., positively homogeneous sets) and locally minimal\footnote{Reducing to the notion of local minimality, rather than absolute minimality, for these cones is necessary, because the (non)local capillarity energy is infinite for cones.}
(i.e., they minimize under compact perturbations, and without volume constraint, the nonlocal capillarity functional, with the original domain~$\Omega$ replaced by its tangent halfspace): see Theorem A.2 in~\cite{MR3717439} and
Corollary 1.3 and~\cite{MR4404780} for precise statements.

Classifying these cones would be of great importance to understand the contact angles
of the nonlocal capillarity minimizers: this is difficult since, in principle, different
subsequences in the above blow-up procedure may produce different limit cones.

In the planar case, however, there is only one possible fractional minimizing cone:
indeed (see Theorem~1.4 in~\cite{MR4404780}) if~$n=2$ and~$E$ is a locally minimizing cone in a halfplane for the nonlocal capillarity functional with
interaction kernel~$K$ is as in~\eqref{KA2s},
then~$E$ is made of only one component,
and there is only one possible contact angle.

In general, the determination of the nonlocal contact angle
requires a blow-up procedure. When the uniqueness of this blow-up is not guaranteed,
we focus at a regular boundary point (say, the origin). In this situation,
we can obtain a nonlocal prescription of the contact angle.
More precisely (according to Theorem~1.4~\cite{MR3717439}),
if~$E$ is a critical set for the nonlocal capillarity functional
with interaction kernel~$K$ is as in~\eqref{KA2s}, up to a translation we can suppose that the origin
lies on~$\partial E$ and we can consider the halfspace~$H$ (respectively,~$V$) such that the blow-up of~$\Omega$ approaches~$H$
(respectively, the blow-up of~$E$ approaches~$V$).
Then, denoting by~$\vartheta$ the angle between~$H$ and~$V$, one finds that
\begin{equation}\label{j1aw} {\mathcal{H}}^{K_s}_{ H\cap V}(v)+(\sigma-1)\int_{\R^n\setminus H}K(v-y)\,dy=0,\end{equation}
for every~$v\in(\partial V)\cap H$. 

Remarkably, this equation uniquely identifies
the angle~$\vartheta=\vartheta(s,\sigma)$ between~$H$ and~$V$, which can thus be seen as a nonlocal contact angle.

Some remarks are in order. First of all, the unique angle for the locally minimizing cones in the plane necessarily satisfies~\eqref{j1aw}.

Furthermore, condition~\eqref{j1aw} looks (and, in fact, is) more complicated than its
classical counterpart (namely, Young's Law in~\eqref{YOU}). Yet, in spite of its involuted formulation, it does determine one and only one contact angle, therefore
it is, somehow unavoidably, the ``only possible'' nonlocal Young's Law.

Besides, as it happened for the Euler-Lagrange equation in~\eqref{el}, the contact angle prescription in~\eqref{j1aw} combines the geometric feature provided by the nonlocal mean curvature with the interaction kernel and the relative adhesion coefficient~$\sigma$.
Actually, philosophically speaking, while in the classical case the Euler-Lagrange equation
(i.e., the constant mean curvature prescription) and Young's Law in~\eqref{YOU} are structurally different, their nonlocal counterparts in~\eqref{el} and~\eqref{j1aw} 
share the very same feature (we may think that the second is just a ``blow-up version'' of the first): this is conceptually interesting, since it highlights the natural ``unifying'' properties of the nonlocal problems.

We also have that, for a given~$s\in(0,1)$, the contact angle is increasing with respect
to the relative adhesion coefficient~$\sigma$, with~$\vartheta(s,-1)=0$, $\vartheta(s,0)=\frac\pi2$ and~$\vartheta(s,1)=\pi$.
This gives that also in the nonlocal case the natural range for the
relative adhesion coefficient is the interval~$[-1,1]$, that~$\sigma>0$ corresponds to the
hydrophobic case of a contact angle larger than~$90^\circ$ and~$\sigma<0$ corresponds to the hydrophilic case of a contact angle smaller than~$90^\circ$.

\subsection{Asymptotic expansions of the nonlocal contact angle}

Interestingly, condition~\eqref{j1aw} formally becomes trivial
in the limit~$s\nearrow1$ 
(simply saying that a halfspace has zero mean curvature), hence the determination of the angle~$\vartheta(s,\sigma)$
relies on (somewhat higher order) integral cancellations.

Looking at these cancellations it is indeed possible to recover the classical Young's Law in~\eqref{YOU} as a limit case as~$s\nearrow1$.
The limit situation as~$s\searrow0$ is also of interest. Interestingly these asymptotics turn
out to be rather explicit and are decribed in Theorem~1.1 of~\cite{MR3707346}. Indeed,
as~$s\nearrow1$ one has that
\begin{equation}\label{as1}
\begin{split}
\vartheta(s,\sigma)\,=\,{
\arccos(-\sigma)-
(1-s)\,\zeta_\sigma+o(1-s)}
,\end{split}
\end{equation}
while as~$s\searrow0$ it holds that
\begin{equation}\label{as0}
\begin{split}
\vartheta(s,\sigma)\,=\,{
\frac\pi2(1+\sigma)-
s\,\eta_\sigma+o(s)}
,\end{split}
\end{equation}
where~$\zeta_\sigma$ and~$\eta_\sigma$ are structural constants which can be computed explicitly.

A direct byproduct of~\eqref{as1} is that~$\cos(\vartheta(s,\sigma))=-\sigma+O(1-\sigma)$
as~$s\nearrow1$
and accordingly~\eqref{j1aw} recovers the classical Young's Law in~\eqref{YOU} as~$s\nearrow1$.

\begin{figure}[h]
    \centering
    \includegraphics[width=6cm]{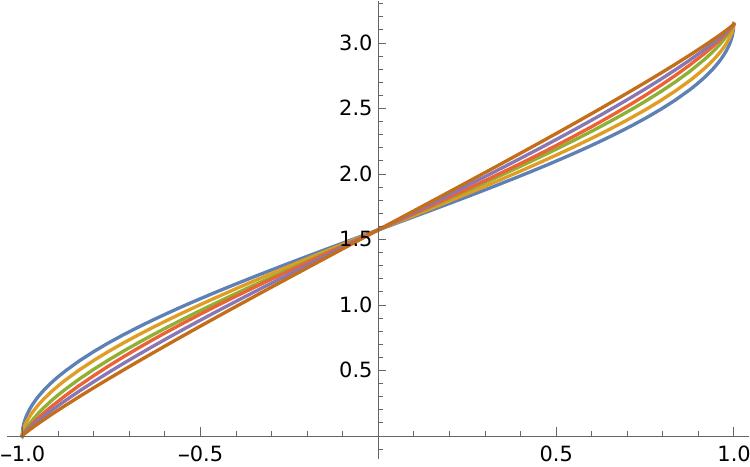}
    \caption{How the dependence of~$\vartheta(s,\sigma)$ upon the
relative adhesion coefficient~$\sigma$ ``linearizes'' when the fractional parameter~$s$
varies from~$1$ to~$0$ (obtained with Mathematica disregarding the higher orders and plotting~\eqref{as1}
with~$s\in\{0.99,\,0.69,\,0.39\}$ and~\eqref{as0}
with~$s\in\{0.31,\,0.21,\,0.11\}$).}
    \label{GOC179}
    \end{figure}

The asymptotics in~\eqref{as0} is also intriguing, because it suggests that the dependence 
of~$\vartheta$ upon~$\sigma$ somewhat ``linearizes'' in the limit as~$s\searrow0$,
since in this case~$\vartheta(s,\sigma)=\frac\pi2(1+\sigma)+O(s)$.
See Figure~\ref{GOC179} for a sketch of these asymptotics.

\subsection{The case of two kernels}
Till now, we have considered the case in which the interaction kernel describes both the interactions between the liquid and the air and those between the liquid and the solid. It is however interesting to consider the case of more general interactions and to allow the model to depend on two different kernels. This scenario has been considered in~\cite{DELUCA}.

Specifically, rather than the setting in~\eqref{ILK}, one can look at the situation in which
two interaction kernels, not necessarily invariant under rotations, are present:
a kernel~$K_1$, accounting for the liquid-air interactions, and a kernel~$K_2$, corresponding to the liquid-solid interactions, such that, for all~$x\in\R^n\setminus\{0\}$ and~$j\in\{1,2\}$,
\begin{eqnarray*}&& \frac{\chi_{B_{\varrho_0}}{(x)}}{\Lambda\,|x|^{n+s_j}}\le K_j(x)\le
\frac{\Lambda}{|x|^{n+s_j}}\\
{\mbox{and}}\qquad&&|\nabla
K_j(x)|\le\frac{\Lambda}{|x|^{n+s_j+1}},\end{eqnarray*}
for some given~$s_1$, $s_2\in(0,1)$.

These kernels give rise to the interactions
$$ {\mathcal{I}}_j(X,Y):={\iint_{X\times Y}K_j(x-y)\,dx\,dy },\qquad j\in\{1,2\}$$
and the corresponding energy functional
\begin{equation}\label{carc}\begin{split}&
{\mathcal{E}}_{\sigma,K_1,K_2}(E):=
{\mathcal{I}}_1(E,\Omega\setminus E)\\&\qquad\qquad\qquad\qquad+\sigma\,
{\mathcal{I}}_2(E,\R^n\setminus\Omega).\end{split}\end{equation}
External forces may be accounted too, but, once again, we focus on the simplest possible case.
As we will see, the situation here is different than before, due to the influence of different scales and the possible lack of perfect cancellations.

In this setting, the Euler-Lagrange equation for a critical set~$E$ of the nonlocal capillarity functional~${\mathcal{E}}_{\sigma,K_1,K_2}$ in~\eqref{carc}
becomes
\begin{equation}\label{el2} \begin{split}&{{\mathcal{H}}^{K_1}_{ E}(x)-\int_{\R^n\setminus\Omega}K_1(x-y)\,dy}\\&\qquad
+\sigma
\int_{\R^n\setminus\Omega}K_2(x-y)\,dy={\rm const},\end{split}\end{equation}see Proposition~1.2 in~\cite{DELUCA}.

Of course, \eqref{el2} boils down to~\eqref{el} when~$K_1=K_2$, but it also clarifies
the role of the coefficient~$(\sigma-1)$ in~\eqref{el}, showcasing that this coefficient indeed arises from a term due to the liquid-air interaction, with a minus sign, plus~$\sigma$ times the contribution due to the liquid-solid interactions (instead, the geometric contribution described by the nonlocal mean curvature is produced only by the liquid-air interaction).\medskip

To understand the contact angle in this case, one considers the blow-up limit of the kernels~$K_1$ and~$K_2$, assuming that, for~$j\in\{1,2\}$, the following limit exists
$$ K^*_j(x):=\lim_{r\searrow0}r^{n+s_j}K_j(x)=\frac{a_j(\vec{x})}{|x|^{n+s_j}},$$
where $j\in\{1,2\},$
for some positive, continuous, even function~$a_j$, where~$\vec{x}:=\frac{x}{|x|}$.

The determination of the contact angle in this case heavily depends on the different homogeneity powers~$s_1$ and~$s_2$
of the kernels
(except in the special case when~$\sigma=0$, in which~$s_2$ does not play any effective role in the minimization procedure).

To appreciate these structural differences, we denote by~$\vartheta(s_1,s_2,\sigma)$ the contact angle in this situation
(i.e., the angle between the halfspace~$H$ tangent to the container~$\Omega$ and the halfspace~$V$ tangent to the droplet~$E$ at the contact point), and (see Theorem 1.4(1)-(2) in~\cite{DELUCA}) we have that, when~${s_1<s_2}$,
\begin{itemize}
\item if~${\sigma<0}$, we have that~${\vartheta(s_1,s_2,\sigma)=0}$,
\item if~${\sigma>0}$, we have that~${\vartheta(s_1,s_2,\sigma)=\pi}$.
\end{itemize}
In a nutshell, when~$s_1<s_2$, the kernel~$K_2$ ``dominates at small scales'' (because is ``more singular'' than the other one): hence, if we believe that
the ``local scale'' produces a prominent effect in the interface as a surface tension,
we expect~$K_2$ to be decisive for the
determination of the contact angle and the kernel~$K_1$ to become ``ineffective''.
Thus, the role of the container becomes decisive and the system basically only sees the interaction of the liquid with the solid.

The situation changes when~$s_1>s_2$ (according to Theorem 1.4(3) in~\cite{DELUCA}). Namely,
if~${s_1>s_2}$ (or if~$s_1\le s_2$ and~$\sigma=0$),
then~${\vartheta(s_1,s_2,\sigma)\in(0,\pi)}$.
Also,
$$ {{\mathcal{H}}^{K_1^*}_{H\cap V}(v)-\int_{\R^n\setminus H}K_1^*(v-y)\,dy=0},$$
for every~$v\in(\partial V)\cap H$ and if
additionally~$a_1$ is constant, then~$\vartheta(s_1,s_2,\sigma)=\frac\pi2$.

The gist here is that
when~$s_1>s_2$, the kernel~$K_1$ ``dominates at small scales'', hence we expect it to be decisive for the
determination of the contact angle (the kernel~$K_2$ becoming ``ineffective'', as it happens when~$\sigma=0$).  
The container here becomes marginal
hence the contact angle is fully determined by the interaction of the liquid with the air.

The more interesting case is thus when~$s_1=s_2$, since the two kernels have a perfect scaling balance
and one expects that both play a role in the determination of the contact angle. 
In this case, we define, for every~$x=(x_1,x_2)\in\partial B_1\subseteq\R^2$
and~$j\in\{1,2\}$,
\begin{equation*}{\footnotesize
a_j^\star(x):=\begin{cases}
a_j(x) \\ \qquad {\mbox{ if }}n=2,\\
\\
\displaystyle
\int_{\R^{n-2}}\frac{a_j\Big(\overrightarrow{
x_1 \,e_1+
x_2 \,e_n+|x|(0,\bar{y},0)}\Big)}{
\big(1+|\bar{y}|^2\big)^{\frac{n+s_j}2}
}\,d\bar{y}
\\ \qquad {\mbox{ if }}n\ge3
\end{cases}}
\end{equation*}
and (see Theorem~1.6, Proposition~1.9 and Theorem 1.10 in~\cite{DELUCA}) we have that
when~${s_1=s_2}$ and
\begin{equation}\label{0oujfn-29roh-32eirj-9034o5t-PK}
|\sigma|<\frac{\displaystyle\int_0^{\pi}a_1^\star(\cos\alpha,\sin\alpha)\,(\sin \alpha)^{s_1}\,d\alpha}{\displaystyle\int_0^{\pi}a_2^\star(\cos\alpha,\sin\alpha)\,(\sin \alpha)^{s_1}\,d\alpha}
\,,\end{equation}
then, ${\vartheta(s_1,s_2,\sigma)\in(0,\pi)}$ and, for every~$v\in(\partial V)\cap H$, the contact angle is characterized by the relation
\begin{equation}\label{GFG}\begin{split}& {{\mathcal{H}}^{K_1^*}_{ H\cap V}(v)-\int_{\R^n\setminus H}K_1^*(v-y)\,dy}\\&\qquad+\sigma
\int_{\R^n\setminus H}K_2^*(v-y)\,dy=0.\end{split}\end{equation}
The bottom line is here that, since the two kernels have similar scaling properties, both participate to the formation of the contact angle, as given by the balanced equation in~\eqref{GFG}. Also, condition~\eqref{0oujfn-29roh-32eirj-9034o5t-PK} is somewhat natural,
since in case of equal kernels reduces to~$|\sigma|<1$, which recovers the classical
condition on the relative adhesion coefficient for the local capillarity theory.
Furthermore, the role of the relative adhesion coefficient showcased in~\eqref{0oujfn-29roh-32eirj-9034o5t-PK} is necessary
to ensure a nontrivial contact angle (as pointed out in Theorems~1.7 and~1.8 in~\cite{DELUCA}).

\section{Conclusion}

Capillarity theory is a beautiful topic in which mathematics naturally meets the applied sciences and in
which technological innovations rely on a deep mathematical understanding of picturesque natural phenomena.

The classical capillarity theory describes the formation of droplets in terms of a surface tension obtained as an average of microscopical forces of adhesion and cohesion.

It is desirable to fully develop a capillarity theory that accounts for these long-range interactions rather than reducing them to their local approximation.

As a first step in this direction, we have introduced and studied a nonlocal
capillarity energy functional modeled on interaction kernels with a power-like decay.
 
Given the complexity of the phenomena involved, the variety of methodologies involved,
and the allure of the natural patterns involved in this study, it is easy to predict that
the interest in these topics will keep increasing in the years to come.

\section*{Acknowledgments}

This work has been supported by the
Australian Laureate Fellowship FL190100081 ``Minimal
surfaces, free boundaries and partial differential equations''.

\end{document}